\documentclass[11pt]{article}
\usepackage[T1]{fontenc}
\usepackage{euscript}
\usepackage{amsmath}
\usepackage{amssymb}
\usepackage{latexsym}
\usepackage{fancyhdr}
\usepackage{graphicx}
\usepackage{graphics}
\usepackage{epstopdf}
\usepackage{epsfig}
\usepackage{diagrams}
\include{amsfonts}

\newarrow{equal}=====
\newarrow{to}---->
\newarrow{onto}----{>>}
\newarrow{embed}>--->
\newarrow{TeXonto}{-}{-}{-}{-}{->>}
\newarrow{TeXto}----{->}
\newarrow{TeXembed}{>}---{->}
\newarrow{TeXembedd}{>}---{->>}

\newcounter{llistadepth}

\newenvironment{manlist}[1]{\addtocounter{llistadepth}{1}
      \edef\llistacontador{llista\romannumeral\the\value{llistadepth}}
      \list{{\rm ({#1{\llistacontador}})}}{\usecounter{\llistacontador}
      \def\makelabel##1{\hss\llap{##1}}
      \itemsep=2pt\parsep=0pt\topsep=3pt plus 1pt minus 1 pt}}{\endlist
      \addtocounter{llistadepth}{-1}}

\renewenvironment{enumerate}{\begin{manlist}{\roman}}{\end{manlist}}

\newtheorem{df}{Definition}[section]      
\newtheorem{thm}[df]{Theorem}             
\newtheorem{prop}[df]{Proposition}
\newtheorem{cor}[df]{Corollary}

\newcommand{\pf}{\noindent{\sc Proof.}\ }

\newcommand{\boom}{\quad\lower3pt\hbox{\vrule height1.1ex width .9ex depth -.2ex}
                    \vskip9pt}

\renewcommand{\mathcal}[1]{\EuScript{#1}}

\let\phi=\varphi

\newcommand{\id}{\mathop{\rm id}\nolimits}
\newcommand{\ad}{\mathop{\rm ad}\nolimits}

\def\Aut{\mathop{\rm Aut}}

\def\id{{\rm id}} 

\def\chigh{{\raise1.5pt\hbox{$\chi$}}}

\newcommand{\gog}{\mathfrak{g}}
\newcommand{\hoh}{\mathfrak{h}}

\newfont{\pointwise}{lcircle10 scaled 500}
\newcommand{\pwise}{{\mbox {\pointwise \symbol{118}}}}
\def\ptwise{\hskip.02in\raise2pt\hbox{$\pwise$}}

\newfont{\numbers}{bbm9 scaled 1200}
\newcommand{\reals}{{\mbox {\numbers R}}}

\newcommand{\co}{\colon\thinspace}                

\def\act{\mathbin{\hbox{$<\kern-.4em\mapstochar\kern.4em$}}}
\def\ract{\mathbin{\hbox{$\mapstochar\kern-.3em>$}}}

\let\Bar=\overline
\let\Hat=\widehat
\let\Tilde=\widetilde

\def\gpd{\,\lower1pt\hbox{$\longrightarrow$}\hskip-.24in\raise2pt
             \hbox{$\longrightarrow$}\,}

\newcommand{\surj}{-\!\!\!-\!\!\!-\!\!\!\gg}
\newcommand{\inj}{>\!\!\!-\!\!\!-\!\!\!-\!\!\!>}

\begin{document}
\title{{\bf Classification of extensions of principal bundles and transitive Lie groupoids with prescribed kernel and cokernel}
\thanks{2000 {\em Mathematics
Subject Classification.} Primary 55R15.
Secondary 22A22, 55N91.
}
\thanks{Keywords: groupoid, extension, cohomology, classification}
\thanks{Research supported by the Greek State Scholarships Foundation
}
\thanks{EDGE Postdoctoral Fellowship, contract number HPRN-CT-2000-00101
}}

\author{Iakovos Androulidakis\\
        Departamento di Mathem\'{a}tica\\
        Instituto Superior T\'{e}cnico\\
        Lisbon 1049-001\\
        Portugal\\
        {\sf iakovos@math.ist.utl.pt}}

\date{{\sf \today}}

\maketitle

\begin{abstract}
The equivalence of principal bundles with transitive Lie groupoids due to Ehresmamm is a well known result. A remarkable generalisation of this equivalence, given by Mackenzie, is the equivalence of principal bundle extensions with those transitive Lie groupoids over the total space of a principal bundle, which also admit an action of the structure group by automorphisms. This paper proves the existence of suitably equivariant transition functions for such groupoids, generalising consequently the classification of principal bundles by means of their transition functions, to extensions of principal bundles by an equivariant form of \v{C}ech cohomology.
\end{abstract}

\section*{Introduction}

Lie groupoids are categories where every arrow has an inverse, plus a smooth structure. They generalise at the same time the notion of a manifold and a group, 
and are widely understood to be part of the general context of 
noncommutative geometry. First, because groupoids are inherently noncommutative 
objects, to a greater extent than are groups. Secondly, Lie groupoids provide 
a modern context for the understanding of the geometry of symplectic and 
Poisson manifolds, which are equipped with noncommutative structures. 
Following a result of Mackenzie, it was shown in \cite{Panh}, that 
the prequantization problem for a symplectic manifold amounts to the existence 
of a suitable transitive Lie groupoid. Furthermore, given a Poisson manifold, 
the existence of a (non-transitive) symplectic groupoid provides a way to 
quantize it.

A rough and descriptive definition of a Lie groupoid is a pair of manifolds $\Omega$ and $M$ such that the elements of $\Omega$ are arrows between points of $M$. The functions $\alpha, \beta \co \Omega \rightarrow M$ mapping every arrow to its source and target points in $M$ are differentiable. Moreover there is a differentiable way to multiply suitable arrows (such that the source of one is exactly the target of the other), and the inversion of arrows is also differentiable. In this setting, for $x, y \in M$ we denote $\Omega_{x}$ the set of arrows in $\Omega$ with source $x$, $\Omega^{y}$ the arrows with target $y$ and $\Omega_{x}^{y}$ the arrows with source $x$ and target $y$. In particular, $\Omega_{x}^{x}$ is a Lie group called the {\em orbit} of $\Omega$ at $x$. A Lie groupoid is denoted by $\Omega \gpd M$.

The simplest example of a Lie groupoid is the product $M \times M \gpd M$ of a manifold $M$, with the obvious groupoid structure. This is called the "pair" groupoid. If $\Omega$ and $\Xi$ are Lie groupoids over the same base manifold $M$, then a smooth map $\phi : \Omega \rightarrow \Xi$ is a {\em morphism} of Lie groupoids if $\alpha \circ \phi = \alpha$, $\beta \circ \phi = \beta$ and $\phi(\eta \cdot \xi) = \phi(\eta) \cdot \phi(\xi)$ for any pair of composable arrows in $\Omega$. For example, given any Lie groupoid $\Omega \gpd M$, the map $(\beta, \alpha) \co \Omega \rightarrow M \times M$ is a morphism of Lie groupoids. This particular morphism is called the {\em anchor}. 

The most well-known classification of Lie groupoids is the one of the transitive case. Transitive Lie groupoids are the ones whose anchor is a surjective submersion, in other words there is an arrow between any two points in $M$. The choice of a basepoint $x \in M$ for a transitive Lie groupoid $\Omega \gpd M$ gives rise to the principal bundle $\Omega_{x}(M, \Omega_{x}^{x}, \beta_{x})$. The principal bundles arising from different choices of elements in $M$ are isomorphic. Given a principal bundle $P(M,G,\pi)$ on the other hand, the associated transitive Lie groupoid is the quotient  $\frac{P \times P}{G} \gpd M$. The groupoid structure here is as follows: For an element $\langle u_{2},u_{1} \rangle$, the source is $\pi(u_{1})$ and the target $\pi(u_{2})$. Suitable arrows $\langle u_{2},u_{1} \rangle$ and $\langle u'_{2},u'_{1} \rangle$ such that there exists a $g \in G$ with $u_{1} = u'_{2}g$ can be multiplied by
\[
\langle u_{2},u_{1} \rangle \langle u'_{2},u'_{1} \rangle = \langle u_{2},u'_{1}g \rangle
\]
The inverse of $\langle u_{2},u_{1} \rangle$ is $\langle u_{1},u_{2} \rangle$ and the unit element over an $x \in M$ is $\langle u,u \rangle$ for any $u \in P$ such that $\pi(u) = x$. It is shown in \cite[II\S1]{LGLADG}, that the two processes are mutually inverse.

So transitive Lie groupoids are classified by the well known classification of principal bundles by \v{C}ech cohomology.

A different classification of the transitive case was given by Mackenzie in \cite{Mac:class}. It was shown that if we shift the point of view from the prescription of $\Omega_{x}^{x}$ (for any given basepoint) to the prescription of the Lie group bundle $I\Omega$ over $M$, of orbits, then transitive Lie groupoids are classified by \v{C}ech cohomology with abelian coefficients. This classification is always possible to calculate in contrast with the often non-abelian classification of principal bundles. To achieve this classification, a transitive Lie groupoid is considered as an extension
\[
I\Omega \inj \Omega \stackrel{(\beta,\alpha)}{\surj} M \times M
\]
of the Lie groupoid $M \times M \gpd M$ (with the obvious groupoid structure) by the Lie group bundle $I\Omega$, instead of the principal bundle $\Omega_{x}(M,\Omega_{x}^{x},\beta_{x})$. For example, the groupoid extension associated to a principal bundle $P(M,G)$ is
\begin{eqnarray}\label{extns0}
\frac{P \times G}{G} \inj \frac{P \times P}{G} \surj M \times M,
\end{eqnarray}
where $\frac{P \times G}{G} \rightarrow M$ is the well known gauge group bundle of $P(M,G)$ (where the $G$-action on itself implied is the adjoint). The usual classification of principal bundles by $\check{H}^{1}(M,G)$ is the answer to the problem "given a Lie group $G$ and a manifold $M$, classify all principal bundles $P(M,G)$". Mackenzie's results imply that if we shift the problem to "given a Lie group bundle $F \rightarrow M$ classify all groupoid extensions of $M \times M$ by this bundle", then we get a classification by \v{C}ech cohomology with coefficients in an abelian group which is always computable, instead of $\check{H}^{1}(M,G)$.

Another classification appeared recently by Moerdijk. In \cite{Moerdijk:classification} regular Lie groupoids are classified, i.e. those ones whose orbits have a constant dimension. Many Lie groupoids are regular, for example those arising from regular Poisson manifolds; moreover all transitive Lie groupoids are regular. Extensions appear in this classification as well. Namely, it is shown that regular Lie groupoids are extensions of foliation groupoids by bundles of connected Lie groups, and they are classified as such. In the case of transitive Lie groupoids, the results in \cite{Mac:class} are a variation of the results of Moerdijk in \cite{Moerdijk:classification}. 

The main result of the present paper is the classification of extensions of transitive Lie groupoids by bundles of Lie groups. Denote such an extension
\begin{eqnarray}\label{extns1}
F \inj \Omega \surj \Xi
\end{eqnarray}
where $F$ is a bundle of Lie groups and $\Omega, \Xi$ are Lie groupoids, all of them over the same connected manifold $M$. Due to the equivalence of transitive Lie groupoids with principal bundles, such extensions are equivalent to extensions of principal bundles 
\begin{eqnarray}\label{extns2}
N \inj Q(M,H) \surj P(M,G).
\end{eqnarray}
Here $N$ is a Lie group and the notation implies the existence of an extension of Lie groups 
\[
N \inj H \surj G.
\]
On the other hand, an extension of principal bundles (\ref{extns2}), gives rise to the extension of transitive Lie groupoids over $M$
\[
\frac{Q \times N}{H} \inj \frac{Q \times Q}{H} \surj \frac{P \times P}{G}
\]
Here the quotient $\frac{Q \times N}{H} \rightarrow M$ is the bundle of Lie groups associated to the principal bundle $Q(M,H)$ through the action of $H$ on $N$ by (the restrictions of) inner automorphisms. It is shown in \cite{Mac3} that the two processes are mutually inverse.

From this point of view, the importance of such a classification is more than the generalisation of the classification of transitive Lie groupoids to extensions. The central problem it deals with is the classification of the covering bundles of a given principal bundle $P(M,G)$ with connected base manifold $M$. Less abstract uses of such a classification arise from an abundance of paradigms of extensions of principal bundles (see for example \cite{Mackenzie:onextnsofpbs}). 

The classification of extensions (\ref{extns1}) is made possible using a result of Mackenzie \cite{Mac3}. It was proved that such extensions are equivalent to a special kind of transitive Lie groupoids, the so-called PBG-groupoids. These are transitive Lie groupoids over the total space of a principal bundle which admit an action of the Lie group of the bundle by Lie groupoid isomorphisms. A description of this equivalence is given in Section 1 of this paper. Roughly speaking, the PBG-groupoid that corresponds to (\ref{extns1}) is a Lie groupoid over the principal bundle $P(M,G)$, together with a $G$-action by (Lie groupoid) automorphisms. Thinking in terms of the extension of principal bundles ({\ref{extns2}) corresponding to (\ref{extns1}), this is a remarkable result; because although the Lie group $G$ does not always act on the kernel $N$ (unless $N$ is abelian), due to Mackenzie's result there always exists a Lie groupoid which admits an action of $G$. 

Once this result is well understood, the problem shifts to the classification of PBG-groupoids. The classification we give here is similar to the one given for general transitive Lie groupoids. In that case, the equivalence with principal bundles ensures the existence of transition functions for Lie groupoids, which suffice to classify them by the usual \v{C}ech cohomology. In the case of PBG-groupoids though, it is necessary to encode the group action as well, and the existence of transition functions which keep track of the action is not established.

In this paper it is shown that there exist transition functions for PBG-groupoids which are equivariant in a certain sense. This is a non-standard notion of equivariance which we call {\em isometablicity}. In turn, a non-standard form of equivariance in \v{C}ech cohomology arises. The first isometablic \v{C}ech cohomology then classifies PBG-groupoids. 

Furthermore, a rather old problem is answered. Lie algebroids are the infinitesimal objects that arise from Lie groupoids, remotely related to them like Lie algebras are related to Lie groups. Mackenzie in \cite{LGLADG} gave a classification of transitive Lie algebroids, but it is not clear how this classification integrates to the groupoid level. A reformulation of the isometablic transition functions is given here, which clearly differentiates to the equivariant analogue of the classification given in \cite{LGLADG}.

This paper is structured in the following way: Section 1 is an account of PBG-groupoids and their relation with extensions of Lie groupoids and principal bundles. In Section 2 the relevant connection theory is described, emphasizing on the material that is of use for the scope of this paper. In Section 3 we prove the existence of transition functions which keep track of the group action, and clarify the notion of isometablicity. Section 4 gives the classification of PBG-Lie group bundles. A remarkable result yielding from this is that the local $G$-actions which give rise to the notion of isometablicity are local expressions of the action of $G$ on the Lie group bundle $I\Omega$ of a given PBG-groupoid $\Omega \gpd P(M,G)$. Section 5 contains the proof of the fact that isometablic transition functions indeed classify PBG-groupoids. Section 6 provides the reformulation of iosmetablic transition functions to a form that differentiates to the equivariant analogue of the classification of Lie algebroids given in \cite{LGLADG}. Finally, the formulation of the suitable cohomology groups where the cocycles of isometablic transition functions live is given in Section 7.

\section*{Acknowledgments}

I would like to thank Kirill Mackenzie for all the support and the useful discussions we had when I was a postgraduate student in Sheffield. Also thanks to John Rawnsley and Rui Fernandes for their valuable comments.

\section{Lie groupoid extensions and PBG-groupoids}

This section recalls in short the material from \cite{Mac3} on the correspondence of extensions of transitive Lie groupoids to PBG-groupoids.
\begin{df}\label{df:PBGgpd}
A {\em PBG-groupoid} is a Lie groupoid $\Omega \gpd P$ whose base is the total
space of a principal bundle $P(M,G)$ together with a right action of $G$ on 
the manifold $\Omega$ such that for all $(\xi,\eta) \in \Omega * \Omega$ and 
$g \in G$ we have:
\begin{enumerate}
\item $\beta(\xi \cdot g) = \beta(\xi) \cdot g$ and $\alpha(\xi \cdot g) =
\alpha(\xi) \cdot g$
\item $1_{u \cdot g} = 1_{u} \cdot g$
\item $(\xi \eta) \cdot g = (\xi \cdot g)(\eta \cdot g)$
\item $(\xi \cdot g)^{-1} = \xi^{-1} \cdot g$
\end{enumerate}
\end{df}
The notation $\Omega * \Omega$ stands for the pairs $(\xi,\eta) \in \Omega \times \Omega$ such that $\alpha(\xi) = \beta(\eta)$.
We denote a PBG-groupoid $\Omega$ over the principal bundle $P(M,G)$ by 
$\Omega \gpd P(M,G)$ and the right-translation in $\Omega$ coming from the 
$G$-action by $\Tilde{R}_{g}$ for any $g \in G$. The right-translation in $P$ 
will be denoted by $R_{g}$. The previous definition implies that $\Tilde{R}_{g}$ is an automorphism of the Lie groupoid 
$\Omega$ over the diffeomorphism $R_{g}$ for all $g \in G$. A morphism $\phi$ 
of Lie groupoids between two PBG-groupoids $\Omega$ and $\Omega'$ over the same 
principal bundle is called a morphism of PBG-groupoids if it preserves the 
group actions, namely if $\phi \circ \Tilde{R}_{g} = \Tilde{R}_{g}' \circ \phi$ 
for all $g \in G$. In the same fashion, a PBG-Lie group bundle (PBG-LGB) is a 
Lie group bundle $F$ over the total space $P$ of a principal bundle $P(M,G)$ 
such that the group $G$ acts on $F$ by Lie group bundle automorphisms. We 
denote a PBG-LGB by $F \rightarrow P(M,G)$. It is easy to see that the gauge 
Lie group bundle $I\Omega \rightarrow P$ associated with a PBG-groupoid $\Omega \gpd P(M,G)$ is a PBG-LGB.

Numerous examples of transitive PBG-groupoids and their corresponding extensions can be found in \cite{Mackenzie:onextnsofpbs}. In \cite{Androulidakis:solo1} non-transitive examples are given as well. Transitive PBG-groupoids are the concern of this paper, due to their equivalence with extensions of transitive Lie groupoids (or, equivalently, extensions of principal bundles  \cite{Mac3}). Let us give an outline of this equivalence.

Given an extension of Lie groupoids (\ref{extns1}), the choice of a basepoint gives rise to its corresponding principal bundle extension (\ref{extns2}) as was discussed in the Introduction. With the notation of (\ref{extns2}), the Lie group $N$ acts on the manifold $Q$ by the restriction of the $H$-action on $Q$ to the embedding of $N$ in $H$. It is immediate that $Q(P,N,\pi)$ is a principal bundle. Here the projection $\pi : Q \surj P$ is the surjective submersion given with the extension (\ref{extns2}). In \cite{Mac3} this was called the {\em transverse bundle}. 

Denote $\Omega$ the (transitive) Lie groupoid $\frac{Q \times Q}{N} \gpd P$ associated to the transverse bundle, and define a right action of the Lie group $G$ on $\Omega$ by
\[
\langle q_{2},q_{1} \rangle g = \langle q_{2}h,q_{1}h \rangle
\]
where $h \in H$ is any element which projects to $g$. It is trivial to see that this action is well defined and makes $\Omega$ a transitive PBG-groupoid over the principal bundle $P(M,G)$.

It is shown in \cite[1.3]{Mac3} that the Lie group bundle $I\Omega \rightarrow P$ of the orbits of $\Omega$ is isomorphic to the pullback bundle $\pi^*(\frac{Q\times N}{H})$. Therefore the PBG-groupoid $\Omega \gpd P(M,G)$ can be presented canonically in the following form:
\[
\pi^*(\frac{Q\times N}{H}) \inj \Omega \surj P \times P.
\]
Here the injection is
\[
(p,\langle q,n \rangle) \mapsto \langle qnh^{-1},qh^{-1} \rangle
\]
where the element $h \in H$ is chosen so that $\pi(q) = p\pi(h)$. Moreover, it is shown in \cite[1.6]{Mac3}, that $I\Omega$ is a PBG-Lie group bundle over $P(M,G)$, the action of $G$ defined as
\[
(p,\langle q,n \rangle)g = (pg,\langle q,n \rangle).
\]

Conversely, consider given a transitive PBG-groupoid $\Upsilon \gpd P(M,G)$. It follows easily from (i) of \ref{df:PBGgpd} that the action of $G$ is free. In \cite[2.2]{Mac3} it is shown that the criterion of Godement (see \cite[16.10.3]{Dieudonne}) applies, therefore the quotient manifold $\frac{\Upsilon}{G}$ exists and the projection $\sharp : \Upsilon \rightarrow \frac{\Upsilon}{G}$ is a surjective submersion. 

This manifold has a natural Lie groupoid structure with base $M$ defined as follows: Since the source and target projections of $\Upsilon$ are $G$-equivariant, they induce maps $\alpha', \beta' : \frac{\Upsilon}{G} \rightarrow M$, which are surjective submersions because the projection $\sharp$, the projection of the principal bundle $P(M,G)$, as well as the source and target maps of $\Upsilon$ as also. Take $u_{1}, u_{2} \in \Upsilon$ such that $\alpha' (\langle u_{1} \rangle) = \beta'(\langle u_{2} \rangle)$. Then there exists $g \in G$ such that $\alpha(u_{1}) = \beta(u_{2})g$, so it is meaningful to define
\[
\langle u_{1} \rangle\langle u_{2} \rangle = \langle u_{1}u_{2}g \rangle.
\]

Finally, the map $(\beta,\alpha) : \Upsilon \rightarrow P \times P$ is equivariant, so it induces a smooth submersion $\pi : \frac{\Upsilon}{G} \rightarrow \frac{P \times P}{G}$. It is clear that this is a groupoid morphism over $M$, and its kernel is $\frac{I\Upsilon}{G}$. Therefore
\[
\frac{I\Upsilon}{G} \inj \frac{\Upsilon}{G} \stackrel{\pi}{\surj} \frac{P \times P}{G}
\]
is an extension of Lie groupoids over $M$. Finally, it is easy to see that the two processes are mutually inverese. In \cite{Mac3} the following theorem is proven:
\begin{thm}
The category of transitive Lie groupoid extensions is equivalent to the category of transitive PBG-groupoids. 
\end{thm}

\section{Connections of PBG-groupoids}

An alternative formulation of the connection theory of principal bundles is by using the Atiyah sequence. Given a principal bundle $P(M,G,p)$, it follows from the fact that the bundle projection $p$ is $G$-invariant, that the vector bundle morphism $Tp : TP \rightarrow TM$ 
quotients to a map $p^* : \frac{TP}{G} \rightarrow TM$ which, like $Tp$, is a fibrewise surjective vector bundle morphism, therefore a surjective submersion. The kernel of this map is of course $\frac{T^{p}P}{G}$, where $T^{p}P$ is the vertical subbundle of $TP$, i.e. the kernel of $Tp$. Now the map $j : \frac{P \times \gog}{G} \rightarrow \frac{T^{p}P}{G}$ induced by 
\[
P \times \gog \rightarrow TP, \  (u,X) \mapsto T_{1}(m_{u})(X)
\]
(where $m_{u} : G \rightarrow P$ is $g \mapsto ug$) is a vector bundle isomorphism (see \cite[Appendix A, 3.2]{LGLADG}). Note that the $G$-action on $\gog$ implied here is the adjoint. Therefore the principal bundle $P(M,G,p)$ gives rise to the extension of vector bundles
\begin{eqnarray}\label{Atiyah}
\frac{P \times \gog}{G} \stackrel{j}{\inj} \frac{TP}{G} \stackrel{p^{*}}{\surj} TM
\end{eqnarray}
which is known as the {\em Atiyah sequence}.

The properties of a connection 1-form $\tilde{\gamma} : TP \rightarrow M \times \gog$ allow it to quotient to a left-splitting $\bar{\gamma} : \frac{TP}{G} \rightarrow \frac{P \times \gog}{G}$ of (\ref{Atiyah}). In turn, the rule
\[
j \circ \bar{\gamma} + \gamma \circ p^{*} = 0
\]
corresponds $\bar{\gamma}$ to a right-splitting $\gamma : TM \rightarrow \frac{TP}{G}$ of (\ref{Atiyah}). This way the connection forms of a principal bundle correspond to the right-splittings of its Atiyah sequence. Respectfully, the curvature of the connection 1-form $\Tilde{\gamma}$ corresponds to the 2-form $R_{\gamma} : TM \times TM \rightarrow \frac{P \times \gog}{G}$ defined by $C_{\gamma}(X,Y) = \gamma[X,Y] - [\gamma(X),\gamma(Y)]$.

The module of sections of the vector bundle $\frac{TP}{G} \rightarrow M$ can be identified with the $G$-invariant vector fields of $P$ (see \cite[Appendix A]{LGLADG}), thus inheriting a Lie bracket which, together with $p^*$, satisfiy the properties of the following definition:
\begin{df}
A {\em Lie algebroid} is a vector bundle $A$ on base $M$ together with a vector bundle map $\sharp : A \rightarrow TM$, called the {\em anchor} of $A$, and a bracket 
$[\ ,\ ] \co \Gamma A \times \Gamma A \rightarrow \Gamma A$ which is $\reals$-bilinear, alternating, satisfies the Jacobi identity, and is such that
\begin{enumerate}
\item $\sharp([X,Y]) = [\sharp X,\sharp Y]$,
\item $[X,fY] = f[X,Y] + (\sharp X)(f)Y$
\end{enumerate}
for all $X,Y \in \Gamma A$ and $f \in C^{\infty}(M)$. 
\end{df}
Basic material on Lie algebroids can be found in \cite{LGLADG} and \cite{Mackenzie:constronalgds}. The notion of a Lie algebroid generalises that of the tangent bundle $TM$ of a given manifold $M$, which can be thought of as a Lie algebroid with the well known Lie bracket of vector fields and the identity as the anchor map. Moreover, any bundle of Lie algebras is a Lie algebroid with zero as the anchor map.

If $A$ and $A'$ are Lie algebroids over the same
base $M$, then a morphism of Lie algebroids $\phi : A \rightarrow A'$ over
$M$ is a vector bundle morphism such that $\sharp' \circ \phi = \sharp$ and 
$\phi([X,Y]) = [\phi(X),\phi(Y)]$ for $X,Y \in \Gamma A$. A Lie algebroid is 
called {\em transitive} if its anchor map is a surjective submersion. In this case the kernel of the anchor map is a bundle of Lie algebras, called the {\em adjoint bundle}, and the Lie algebroid can be presented as an extension of vector bundles
\begin{eqnarray}\label{*}
L \inj A \stackrel{\sharp}{\surj} TM
\end{eqnarray}
where the injection of $L$ into $A$ and the anchor map are morphisms of Lie algebroids. 
\begin{df}
Let $A, A'$ be Lie algebroids over the manifold $M$ and $L \rightarrow M$. An extension of vector bundles
\[
K \inj A \surj A'.
\]
is called an {\em extension of Lie algebroids} if the injection and surjection maps are morphisms of Lie algebroids.
\end{df}
Extensions such as (\ref{*}) are the simplest form of Lie algebroid extensions, in fact they are just an alternative way to present a transitive Lie algebroid $A$ over a manifold $M$. In this setting, the connection theory of principal bundles gives rise to the following notions:
\begin{df}
Let $L \inj A \stackrel{\sharp}{\surj} TM$ be a transitive Lie algebroid. 
\begin{enumerate}
\item A {\em connection} of $A$ is a vector bundle morphism $\gamma : TM \rightarrow A$ such that $\sharp \circ \gamma = 0$.
\item The {\em curvature} of a connection $\gamma$ is the 2-form $C_{\gamma} : TM \times TM \rightarrow L$ defined by
\[
C_{\gamma}(X,Y) = \gamma[X,Y] - [\gamma(X),\gamma(Y)]
\] 
for all $X, Y \in \Gamma A$.
\end{enumerate}
A connection $\gamma$ is called {\em flat} if $C_{\gamma} = 0$.
\end{df}
Note that a flat connection is evidently a morphism of Lie algebroids $\gamma : TM \rightarrow A$.

All Lie groupoids differentiate to Lie algebroids. A full account of this process can be found in \cite[III\S3]{LGLADG}. The reader can get a rough idea by comparing the extension (\ref{extns0}) to the Atiyah sequence (\ref{Atiyah}). Lie III does not apply for groupoids and algebroids though. The integrability of Lie algebroids has a cohomological obstruction in the transitive case, which was given by Mackenzie in \cite[V]{LGLADG}. In the non-transitive case, integrability of Lie algebroids is a problem of different order which was tackled by Crainic and Fernandes in \cite{CF}. In general, a Lie algebroid that integrates to a Lie groupoid $\Xi \gpd M$ is denoted by $A\Xi$. Note that the tangent bundle $TM$ of a manifold $M$ integrates to the "pair" groupoid $M \times M \gpd M$.

Analogously with the reformulation of principal bundle connections as right-splittings of the Atiyah sequence, it is legitimate to regard the  {\em connections} of a transitive Lie groupoid $\Xi \gpd M$ as the connections of the Lie algebroid $A\Xi$ it differentiates to, and the same is valid for the curvature 2-forms. This terminology will be used in the remaining of this paper. 

Once again though, the concern of this paper is extensions of transitive Lie groupoids, so let us make a fresh start by giving the notion of a PBG-algebroid.

\begin{df}\label{df:PBGalgd}
A {\em PBG-algebroid} over the principal bundle $P(M,G)$ is a Lie 
algebroid $A$ over $P$ together with a right action of $G$ on $A$ denoted by 
$(X,g) \mapsto \Hat{R}_{g}(X)$ for all $X \in A,~ g \in G$ such that each 
$\Hat{R}_{g} : A \rightarrow A$ is a Lie algebroid automorphism over the right 
translation $R_{g}$ in $P$. 
\end{df}

We denote a PBG-algebroid $A$ over $P(M,G)$ by $A \Rightarrow P(M,G)$. The $G$-action 
on $A$ induces an action of $G$ on the module $\Gamma A$ of sections of the vector bundle $A \rightarrow M$, namely
\[
X \cdot g = \Hat{R}_{g} \circ X \circ R_{g^{-1}}.
\]
The right-translation with respect to this action is denoted by $\Hat{R}_{g}^{\Gamma} \co \Gamma A \rightarrow \Gamma A$ for all $g \in G$. With this notation definition \ref{df:PBGalgd} implies that 
\[
\Hat{R}_{g}^{\Gamma}([X,Y]) = [\Hat{R}^{\Gamma}_{g}(X),\Hat{R}^{\Gamma}_{g}(Y)].
\] 

Given a {\em transitive} PBG-algebroid $A \Rightarrow P(M,G,p)$, its adjoint bundle $L \rightarrow P$ inherits a $G$-action by automorphisms, thus making 
\[
L \inj A \stackrel{\sharp}{\surj} TP
\]
an extension of PBG-algebroids. That is to say it is an extension of Lie algebroids such that the injection and surjection maps are moreover equivariant. It is shown in \cite[3.4]{Androulidakis:solo1} that the Godement criterion applies, so the quotient manifold $\frac{A}{G}$ exists. Therefore the previous extension quotients to a vector bundle extension
\begin{eqnarray}\label{PBGtoextn}
\frac{L}{G} \inj \frac{A}{G} \stackrel{\sharp^{/G}}{\surj} \frac{TP}{G}
\end{eqnarray}
of the (integrable) Lie algebroid $\frac{TP}{G}$ by the quotient Lie algebra bundle $\frac{L}{G}$. Observe that since the quotient manifold $\frac{A}{G}$ exists, the vector bundle structure of $A$ quotients to $\frac{A}{G} \rightarrow M$. Moreover, the natural projection $\natural^A \co A \rightarrow \frac{A}{G}$ is a pullback over $p \co P \rightarrow M$. 

The vector bundle $\frac{A}{G}$ has the following Lie algebroid structure: The anchor is the composition of vector bundle morphisms $p^* \circ \sharp^{/G}$. Moreover, the sections of $\frac{A}{G}$ are isomorphic to the $G$-invariant sections of $A$, therefore $\Gamma (\frac{A}{G})$ inherits the Lie bracket from $\Gamma^{G}A$. The verification that this bracket together with the anchor map $p^{*} \circ \sharp^G$ satisfy the properties of a Lie algebroid can be found in \cite[3.2]{Mac3}. It is immediate that $\frac{A}{G}$ is transitive. A more elaborate presentation of the extension \ref{PBGtoextn} is given in Figure \ref{fullPBGtoextn}, which helps to keep track of all the structures related to the Lie algebroid extension. Note that the adjoint bundle $K$ of $\frac{A}{G}$ is an extension of $\frac{P \times \gog}{G}$ by $\frac{L}{G}$. This diagram makes clear that the cokernel of the extension (\ref{PBGtoextn}) is in fact the Atiyah sequence of the bundle $P(M,G,p)$.

\begin{figure}[htb]
\begin{eqnarray*}
\begin{diagram}
 &  & K &  & \frac{P \times \gog}{G} \\
 &  & \dTeXembed &  & \dTeXembed>{j} \\
 \frac{L}{G} & \rTeXembed & \frac{A}{G} & \rTeXonto^{\sharp^{/G}} & \frac{TP}{G} \\
 &  & \donto<{p^* \circ \sharp^{/G}} &  & \donto>{p^*} \\
 &  &  TM  & = & TM 
\end{diagram}
\end{eqnarray*}
\caption{The extension of Lie algebroids induced by a PBG-algebroid\label{fullPBGtoextn}}
\end{figure}

On the other hand, pulling back (\ref{PBGtoextn}) by the map $Tp \co TP \rightarrow TM$ we recover the given PBG-algebroid (see \cite[4]{Mac3}). This consists the proof of the following theorem.
\begin{thm}
The category of transitive PBG-algebroids over a manifold $M$ is equivalent to the category of  Lie algebroid extensions
\begin{eqnarray}\label{Lalgextns}
K \inj A \surj A\Xi
\end{eqnarray}
of an integrable transitive Lie algebroid by a Lie algebra bundle (over $M$).
\end{thm}

Now extensions of Lie groupoids differentiate to extensions (\ref{Lalgextns}). The connection theory of Lie groupoid extensions (\ref{extns1}) is encoded by the right-splittings of extensions (\ref{Lalgextns}). These in turn correspond to the following notion of connection for the equivalent PBG-algebroid (see \cite{Mackenzie:onextnsofpbs} and \cite{Mac3}):
\begin{df}\label{conns:PBG}
Let $A \Rightarrow P(M,G,p)$ be a transitive PBG-algebroid. A connection 
$\gamma : TP \rightarrow A$ is called {\em isometablic}, if it satisfies 
\begin{equation}\label{eq:isomet}
\gamma \circ TR_{g} = \Hat{R}_{g} \circ \gamma
\end{equation}
\end{df}

An account of isometablic connections and their holonomy is given in \cite{Androulidakis:solo1}, however in this paper we are interested in a different problem. The groupoid extensions that we intend to classify have prescribed kernel and cokernel. In other words, {\em given} a transitive Lie groupoid $\Xi \gpd M$ and a Lie algebra bundle $F \rightarrow M$, we classify all transitive Lie groupoids $\Phi \rightarrow M$ which fit into a Lie groupoid extension
\[
F \inj \Phi \surj \Xi.
\]
In this sense, we are interested in the connections of $\Phi$ rather than the splittings of the extension of the Lie algebroid extension $AF \inj A\Phi \surj A\Xi$. The following theorem clarifies exactly what these connections correspond to in the relevant PBG-algebroid.
\begin{thm}\label{basisom}
Suppose given a transitive PBG-algebroid $A \Rightarrow P(M,G,p)$ and consider its corresponding extension of Lie algebroids (\ref{PBGtoextn}) over $M$. The connections of the (transitive) Lie algebroid $\frac{A}{G} \rightarrow M$ are equivalent to the isometablic connections of $A$ which vanish on the kernel $T^{p}P$ of $Tp : TP \rightarrow TM$.
\end{thm}

\pf Consider an isometablic connection $\gamma : TP \rightarrow A$ such that $\gamma(X) = 0$ if $X \in T^{p}P$. This quotients to a splitting $\gamma^{/G} : \frac{TP}{G} \rightarrow \frac{A}{G}$. Given a connection $\delta : TM \rightarrow \frac{TP}{G}$ of the principal bundle $P(M,G)$, define
\[
\Tilde{\gamma} = \gamma^{/G} \circ \delta : TM \rightarrow \frac{A}{G}.
\]
The assumption that $\gamma$ vanishes on the kernel of $Tp$ makes the definition of $\Tilde{\gamma}$ independent from the choice of $\delta$. It follows immediately from the assumption that $\delta$ is a connection of $P(M,G)$ and $\gamma^{/G}$ is a splitting of (\ref{PBGtoextn}) that this is a connection of the Lie algebroid $\frac{A}{G}$.

Conversely, given a connection $\theta : TM \rightarrow \frac{A}{G}$ of the Lie algebroid $\frac{A}{G}$, compose it with the anchor map $p^* : \frac{TP}{G} \surj TM$ of the Atiyah sequence corresponding to the bundle $P(M,G,p)$ (see (\ref{Atiyah})) to the vector bundle morphism
\[
\Bar{\theta} = \theta \circ p^* : \frac{TP}{G} \rightarrow \frac{A}{G}.
\]
Denote $\natural : TP \rightarrow \frac{TP}{G}$ and $\natural^A : A \rightarrow \frac{A}{G}$ the natural projections. Since $\natural^A$ is a pullback over $p \co P \rightarrow M$, there is a unique vector bundle morphism $\gamma : TP \rightarrow A$ such that
\[
\natural^A \circ \gamma = \Bar{\theta} \circ \natural.
\]
Due to the $G$-invariance of $\natural$ and $\natural^A$the morphism of vector bundles $\Hat{R}_{g^{-1}} \circ \gamma \circ TR_{g}$ also satisfies the previous equation for every $g \in G$, therefore it follows from the uniqueness argument that $\gamma$ is isometablic. It is an immediate consequence of the previous equation that $\gamma$ vanishes at $T^pP$. 

To see that it is indeed a connection of $A$, let us recall the fact that $\theta$ is a connection of $\frac{A}{G}$. This gives $p^* \circ \sharp^{/G} \circ \theta = \id_{TM}$. Now $\sharp^{/G} = \natural \circ \sharp$ and by definition we have $p^{*} \circ \natural = Tp$, therefore $Tp \circ \sharp \circ \theta = \id_{TM}$. Now take an element $X \in TP$. Then $Tp(X) \in TM$, and it follows from this equation that there exists an element $g \in G$ such that 
\[
(\sharp \circ \theta)(Tp(X)) = X \cdot g.
\]
Multiplying this by $g^{-1}$ and using the $G$-invariance of $Tp$ we get
\[
\sharp \circ (\theta \circ Tp) = \id_{TP}.
\]
Finally, from the properties of the pullback, it follows immediately that $\gamma$ is the map $(\pi, \Bar{\theta} \circ \natural)$, where $\pi \co TP \rightarrow P$ is the natural projection of the tangent bundle. It is straightforward to check that this reformulates to $(\pi, \theta \circ Tp)$, and this proves that $\gamma$ is a connection.
\boom

\begin{df}
The isometablic connections of a PBG-algebroid $A \Rightarrow P(M,G,p)$ which vanish at the kernel $T^{p}P$ of $p^{*}$ are called {\em basic connections}.
\end{df}

It is therefore necessary to focus on basic connections of PBG-groupoids for the purpose of this paper. The following result follows from the proof \ref{basisom}.
\begin{cor}\label{corbasisom}
Let $A \Rightarrow P(M,G)$ be a transitive PBG-algebroid. A flat connection of the Lie algebroid $\frac{A}{G} \rightarrow M$ gives rise to a unique flat basic connection of $A$.
\end{cor}
Note that the proof of \ref{basisom} does not give force to the converse of this result. That is because  the connection of $\frac{A}{G}$ corresponding to a given flat basic connection of $A$ arises by composition with an arbitrary connection of $\frac{TP}{G}$, which is not necessarily a flat one, unless the bundle $P(M,G)$ is flat.

\section{Transition functions for transitive PBG-groupoids}

This section is concerned with the study of those transition functions of transitive PBG-groupoids which encode the group action.

Let us start with a principal bundle $P(M,G)$ and a simple open cover ${\cal{U}} = \{ U_{i} \}_{i \in I}$ of $M$. This is an open cover such that each $U_{i}$ is contractible, and the intersection of two as well as three open sets is also contractible. Then a cover ${\cal{P}} = \{ P_{i} \}_{i \in I}$ of $P$ by principal bundle charts such that $P_{i} \cong U_{i} \times G$ exists. 

Consider now a PBG-groupoid $\Omega \gpd P(M,G)$ over this bundle and its corresponding Lie algebroid $A\Omega \Rightarrow P(M,G)$ with adjoint bundle $L\Omega$. The extension of Lie algebroids corresponding to that is
\[
\frac{L\Omega}{G} \inj \frac{A\Omega}{G} \surj \frac{TP}{G}.
\]
It follows from \cite[IV\S 4]{LGLADG} that the Lie algebroid $\frac{A\Omega}{G}$ (over $M$) has local flat connections $\Tilde{\theta}_{i}^{*} \co TU_{i} \rightarrow (\frac{A\Omega}{G})_{U_{i}}$. Due to \ref{corbasisom} these give rise to flat {\em basic} connections $\theta_{i}^{*} \co TP_{i} \rightarrow A\Omega_{P_{i}}$. 

Since the connections $\Tilde{\theta}_{i}^{*}$ are flat, they can be regarded as morphisms of Lie algebroids. Now consider the following theorem from \cite{Mac-Xu}:
\begin{thm}
Let $\Omega, \Xi$ be Lie groupoids over the same manifold $M$ and $\mu : A\Omega \rightarrow A\Xi$ a Lie algebroid morphism. If $\Omega$ is $\alpha$-simply connected, then there exists a unique morphism of Lie groupoids $\phi : \Omega \rightarrow \Xi$ which differentiates to $\mu$, i.e. $\phi^{*} = \mu$.
\end{thm}
With the assumption that every $U_{i}$ is contractible, and by force of the previous result, it follows that the $\Tilde{\theta}_{i}^{*}$s integrate uniquely to morphisms of Lie groupoids $\Tilde{\theta_{i}} : U_{i} \times U_{i} \rightarrow \frac{\Omega}{G}_{U_{i}}^{U_{i}}$. It was shown in the proof of \ref{basisom} that the basic flat connections $\theta_{i}^{*}$ corresponding to the $\Tilde{\theta}_{i}^{*}$s are in essence the maps $\tilde{\theta}_{i}^{*} \circ Tp$, therefore they also integrate uniquely to morphisms of Lie groupoids
\[
\theta_{i} \co P_{i} \times P_{i} \rightarrow \Omega_{P_{i}}^{P_{i}}.
\]
\begin{prop}
The $\theta_{i}$s are morphisms of PBG-groupoids.
\end{prop}

\pf It suffices to prove the equivariance of the $\theta_{i}$s. For every $g \in G$, the map $\theta_{i}^{g} : P_{i} \times P_{i} \rightarrow \Omega_{P_{i}}^{P_{i}}$ defined by
\[
\theta_{i}^{g}(u,v) = \theta_{i}(ug,vg)g^{-1}
\]
is clearly a morphism of Lie groupoids and it differentiates to $\theta_{i}^{*}$. It therefore follows from the uniqueness of $\theta_{i}$ that $\theta_{i}^{g} = \theta_{i}$ for all $g \in G$, consequently $\theta_{i}$ is equivariant. \boom

For every $i \in I$ choose an element $u_{i} \in P_{i}$ and define $\Bar{\sigma}_{i} \co P_{i} \rightarrow \Omega_{P_{i}}$ by $\Bar{\sigma}_{i}(u) = \theta_{i}(u,u_{i})$. We call these maps {\em schisms}. Note that $\Bar{\sigma}_{i}(u_{i}) = 1_{u_{i}}$. The following proposition clarifies the behaviour of the schisms with respect to the $G$-action. We call this notion of equivariance  {\em isometablicity} because it follows directly from the isometablicity property of the local flat connections of the PBG-groupoid we discussed above.
\begin{prop}
The schisms $\Bar{\sigma}_{i}$ are {\em isometablic} in the sense
\[
\Bar{\sigma}_{i}(ug) = (\Bar{\sigma}_{i}(u)g) \cdot \Bar{\sigma}_{i}(u_{i}g)
\]
for all $u \in P_{i}$ and $g \in G$.
\end{prop}

\pf From the definition of the $\Bar{\sigma}_{i}$s and the equivariance of the morphisms $\theta_{i}$ we get:
\[
(\Bar{\sigma}_{i}(u)g) \cdot \Bar{\sigma}_{i}(u_{i}g) = (\theta_{i}(u,u_{i})g) \cdot \theta_{i}(u_{i}g,u_{i}) = \theta_{i}(ug,u_{i}g) \cdot \theta_{i}(u_{i}g,u_{i}) 
= \theta_{i}(ug,u_{i}) = \Bar{\sigma}_{i}(ug).
\]
\boom

For every choice of a $u_{i} \in P_{i}$, consider the Lie group $H_{i} = \Omega_{u_{i}}^{u_{i}}$. In order to refer to a unique Lie group independent to the index $i \in I$, we need to fix a $u_{0} \in P$ and define $H = \Omega_{u_{0}}^{u_{0}}$. Then, for every $i \in I$ choose a $\xi_{i} \in \Omega_{u_{0}}^{u_{i}}$ and consider the maps $\tau_{i} \co H_{i} \rightarrow H$ defined by $\tau_{i}(\eta) = \xi^{-1} \cdot \eta \cdot \xi$. These are isomorphisms of Lie groups. Now define $\sigma_{i} \co P_{i} \rightarrow \Omega_{u_{0}}$ by 
\[
\sigma = \Bar{\sigma}_{i} \cdot \xi_{i}.
\]
These are sections of the Lie groupoid $\Omega$. Note that $\sigma_{i}(u_{i}) = \xi_{i}$. The isometablicity of these sections is described in the following proposition:
\begin{prop}
The sections $\sigma_{i}$ are {\em isometablic} in the sense
\[
\sigma_{i}(ug) = [\sigma_{i}(u)g] \cdot (\xi_{i}^{-1}g) \cdot \sigma_{i}(u_{i}g).
\]
for all $i \in I, u \in P_{i}$ and $g \in G$.
\end{prop}
The proof is a straightforward calculation. 

Now we look at the isometablicity of the transition functions. We denote $\{ \Bar{s}_{ij} \co P_{ij} \rightarrow \Omega_{u_{j}}^{u_{i}} \}_{i,j \in I}$ the transition functions of the schisms $\{ \Bar{\sigma}_{i} \}_{i \in I}$ and $\{ s_{ij} \co P_{ij} \rightarrow \Omega_{u_{0}}^{u_{0}} \}_{i,j \in I}$ the transition functions of the sections $\{ \sigma_{i} \}_{i \in I}$. The following proposition is an immediate consequence of the isometablicity of the schisms and the sections.
\begin{prop}\label{isom:trans}
For every $i,j \in I$ such that $P_{ij} \neq \emptyset, u \in P_{ij}$ and $g \in G$ we have:
\begin{enumerate}
\item $\Bar{s}_{ij}(ug) = \Bar{\sigma}_{i}(u_{i}g)^{-1} \cdot (\Bar{s}_{ij}(u)g) \cdot \Bar{\sigma}_{j}(u_{j}g)$ 
\item $s_{ij}(ug) = \sigma_{i}(u_{i}g)^{-1} \cdot (\xi_{i}g)\cdot (s_{ij}(u)g) \cdot (\xi_{j}g)^{-1} \cdot \sigma_{j}(u_{j}g)$.
\end{enumerate}
\end{prop}
This gives rise to the following formulation of $G$-actions:
\begin{df}
Denote $\Omega_{u_{i}}^{u_{i}} = H_{i}$ and $\Omega_{u_{0}}^{u_{0}} = H$. The formulas 
\begin{enumerate}
\item $\Bar{\rho}_{ij} \co G \times H_{i} \rightarrow H_{i}$, $\Bar{\rho}_{ij}(g^{-1})(h_{i}) = \Bar{\sigma}_{i}(u_{i}g)^{-1} \cdot (h_{i}g) \cdot \Bar{\sigma}_{j}(u_{j}g)$ and
\item $\rho_{ij} \co G \times H \rightarrow H$, $\rho_{ij}(g^{-1})(h) = \sigma_{i}(u_{i}g)^{-1} \cdot (\xi_{i}g)\cdot (hg) \cdot (\xi_{j}g)^{-1} \cdot \sigma_{j}(u_{j}g)$
\end{enumerate}
define families of $G$-actions on $H_{i}$ and $H$ respectively.
\end{df}
With this notation, it is legitimate to reformulate the isometablicity equations of \ref{isom:trans} to
\[
\Bar{s}_{ij}(ug) = \Bar{\rho}_{ij}(g^{-1})(\Bar{s}_{ij}(u))
\]
and
\[
s_{ij}(ug) = \rho_{ij}(g^{-1})(s_{ij}(u)).
\]

Let us now examine the properties of the $G$-actions $\Bar{\rho}_{ij}$ and $\rho_{ij}$. The proof of the following proposition is, again, straightforward.
\begin{prop}
Let $\Omega \gpd P(M,G)$ be a PBG-groupoid. Then the families of $G$-actions $\{ \Bar{\rho}_{ij} \}_{i,j \in I}$ and $\{ \rho_{ij} \}_{i,j \in I}$ satisfy the following identities:
\begin{eqnarray}\label{cocmorph}
\rho_{ij}(g^{-1})(h_{1}h_{2}) = \rho_{ik}(g^{-1})(h_{1})\rho_{kj}(g^{-1})(h_{2})
\end{eqnarray}
for all $i,j,k \in I$ such that $P_{ijk} \neq \emptyset$ and $h_{1}, h_{2} \in H$.
\begin{multline}\label{indepj}
\rho_{ij}(g^{-1})(h) = \rho_{ii}(g^{-1})(h) \cdot \sigma_{i}(u_{i}g)^{-1} \cdot (\xi_{i}g)\cdot (\xi_{j}g)^{-1} \cdot \sigma_{j}(u_{j}g) = \\
= \sigma_{i}(u_{i}g)^{-1} \cdot (\xi_{i}g)\cdot (\xi_{j}g)^{-1} \cdot \sigma_{j}(u_{j}g) \cdot \rho_{jj}(g^{-1})(h).
\end{multline}
\begin{eqnarray}\label{maptoisomet}
\tau_{i}(\Bar{\rho}_{ii}(g^{-1})(h_{i})) = \rho_{ii}(g^{-1})(\tau_{i}(h_{i}))
\end{eqnarray}
for all $h_{i} \in H_{i}$.
\end{prop}
Due to (\ref{indepj}), it is possible to say that the family of actions $\{\rho_{ij}\}_{i,j \in I}$ is fully determined by the subset of those actions with $i = j$. Now (\ref{maptoisomet}) shows that for all $i \in I$ the isomorphism $\tau_{i} \co H_{i} \rightarrow H$ maps every $G$-action $\Bar{\rho}_{ii}$ on $H_{i}$ exactly to the $G$-action $\rho_{ii}$ on $H$. 

Last, notice that (\ref{cocmorph}) is a non-standard property. From this it follows immediately that $\rho_{ii}(g^{-1})(e_{H}) = e_{H}$ for all $i \in I$. These two properties almost make the $\rho_{ij}$s representations in a certain sense. We single out (\ref{cocmorph}) by giving the following definition.
\begin{df}
Let $G$ and $H$ be Lie groups. If a family $\{ \rho_{ij} \}_{i,j \in I}$ of $G$-actions on $H$ satisfy
\[
\rho_{ij}(g^{-1})(h_{1}h_{2}) = \rho_{ik}(g^{-1})(h_{1})\rho_{kj}(g^{-1})(h_{2})
\]
for all $g \in G$, $h_{1}, h_{2} \in H$ and $i,j,k \in I$ then $G$ is said to be acting on $H$ by {\em cocycle morphisms}.
\end{df}

\subsection{Equivalence of transition functions}

So far we have demonstrated that PBG-groupoids have sections which are suitably equivariant. These sections arise naturally from the local flat basic connections that exist on the algebroid level. But what happens if we start with a different family of local flat basic connections? 

Let us start with two families $\{ \theta_{i}^{*} \}_{i \in I}$ and $\{ {\theta'}_{i}^{*} \}_{i \in I}$ of flat basic connections over the same cover ${\cal{P}} = \{ P_{i} \}_{i \in I}$ of $P$ by principal bundle charts. Then there exist maps $\ell^{*}_{i} \co TP_{i} \rightarrow P_{i} \times \hoh_{i}$ such that 
\[
{\theta'}_{i}^{*} = \theta_{i}^{*} + \ell_{i}^{*}
\]
for every $i \in I$. Here $\gog_{i}$ denotes the Lie algebra of the Lie group $H_{i}$. Therefore every $\ell_{i}^{*}$ must also be isometablic, that is to say
\[
\ell^{*}_{i}(Xg) = \ell^{*}_{i}(X)g
\]
for all $X \in TP_{i}$ and $g \in G$. Moreover, the $\ell_{i}^{*}$s integrate to PBG-groupoid morphisms $\ell_{i} \co P_{i} \times P_{i} \rightarrow H_{i}$ such that ${\theta'}_{i} = \theta_{i} \cdot \ell_{i}$. As far as the isometablicity of the $\ell_{i}$s is concerned, it follows
\begin{eqnarray}\label{isomell}
\ell_{i}(ug,vg) = \Bar{\rho}_{ii}(g^{-1})(\ell_{i}(u,v)).
\end{eqnarray}

Now define $\Bar{r}_{i} \co P_{i} \rightarrow H_{i}$ by 
\[
\Bar{r}_{i}(u) = \ell_{i}(u,u_{i})
\]
and $r_{i} \co P_{i} \rightarrow H$ by $r_{i} = \tau_{i}\circ \Bar{r}_{i}$. That is to say,
\[
r_{i}(u) = \xi_{i}^{-1} \cdot \Bar{r}_{i}(u) \cdot \xi_{i}
\]
for all $u \in P_{i}$. We call the $\Bar{r}_{i}$s and the $r_{i}$s {\em conjugation} maps. The proof of the following proposition is a simple calculation.
\begin{prop}
The schisms, sections and the respective transition data induced by $\{ \theta_{i}^{*} \}_{i \in I}$ and $\{ {\theta'}_{i}^{*} \}_{i \in I}$ are related by:
\begin{enumerate}
\item $\Bar{\sigma'}_{i} = \Bar{\sigma}_{i} \cdot \Bar{r}_{i}$.
\item $\Bar{s'}_{ij} = \Bar{r}_{i}^{-1} \cdot \Bar{s}_{ij} \cdot \Bar{r}_{j}$.
\item ${\sigma'}_{i} = \sigma_{i} \cdot r_{i}$.
\item ${s'}_{ij} = r_{i}^{-1} \cdot s_{ij} \cdot r_{j}$.
\end{enumerate}
\end{prop}
\begin{cor}
The families of $G$-actions ${\cal{\rho}} = \{ \rho_{ij} \}_{i,j \in I}$ and ${\cal{\rho'}} = \{ {\rho'}_{ij} \}_{i,j \in I}$ arising from the connections $\theta_{i}^{*}$ and ${\theta'}_{i}^{*}$ respectively are related by
\[
{\rho'}_{ij}(g^{-1})(h) = r_{i}(u_{i}g)^{-1} \cdot \rho_{ij}(g^{-1})(h) \cdot r_{j}(u_{j}g)
\]
for all $h \in H$ and $g \in G$.
\end{cor}

Now let us examine the isometablicity of the conjugation maps.
\begin{prop}\label{isomconj}
The conjugation maps satisfy
\begin{enumerate}
\item $\Bar{r}_{i}(ug) = \Bar{\rho}_{ii}(g^{-1})(\Bar{r}_{i}(u)) \cdot \Bar{r}_{i}(u_{i}g)$
\item $r_{i}(ug) = \rho_{ii}(g^{-1})(r_{i}(u)) \cdot r_{i}(u_{i}g)$
\end{enumerate}
for all $u \in P_{i}$ and $g \in G$.
\end{prop}

\pf Note that (ii) follows by applying the isomorphisms $\tau_{i}$ to (i) and taking into account (\ref{maptoisomet}). For (i) we have:
\[
\Bar{r}_{i}(ug) = \ell_{i}(ug,u_{i}) = \ell_{i}(ug,u_{i}g) \cdot \ell(u_{i}g,u_{i})
\]
Because of (\ref{isomell}) the last part of the above equation becomes $\Bar{\rho}_{ii}(g^{-1})(\ell(u,u_{i})) \cdot \ell(u_{i}g,u_{i})$, and the result follows. \boom

\section{The classification of PBG-Lie group bundles}
\begin{sloppypar}
Consider the adjoint bundle $I\Omega \rightarrow P(M,G)$ associated with a  a given PBG-groupoid $\Omega \gpd P(M,G)$. This section is concerned with the isometablic transition data that classifies this bundle. Apart from this classification, another result given here is that by using this data it is shown that the $G$-actions $\rho_{ij}$ given in the previous section are local expressions of the action of $G$ on the Lie group bundle $I\Omega$.
\end{sloppypar}
\begin{prop}
Let $\{ U_{i} \}_{i \in I}$ be a simple open cover of $M$ and $P_{i} \cong U_{i} \times G$ charts of the principal bundle $P(M,G)$. The maps $\psi_{i} \co P_{i} \times H \rightarrow I\Omega_{P_{i}}$ defined by
\[
\psi_{i}(u,h) = \sigma_{i}(u) \cdot h \cdot \sigma_{i}(u)^{-1}
\]
are local charts for the Lie group bundle $I\Omega$. They are {\em isometablic} in the sense
\[
\psi_{i}(ug,\rho_{ii}(g^{-1})(h)) = \psi_{i}(u,h) \cdot g.
\]
\end{prop}

\pf The fact that $\psi_{i}$ is a bijection and $\psi_{i,u} \co H \rightarrow \Omega_{u}^{u}$ is a morphism of Lie groups for all $u \in P_{i}$ are simple calculations. For the isometablicity we have:
\begin{multline*}
\psi_{i}(ug,\rho_{ii}(g^{-1})(h)) = \sigma_{i}(ug)\cdot \rho_{ii}(g^{-1})(h) \cdot \sigma_{i}(ug)^{-1} = \\
= (\sigma_{i}(u)g) \cdot (\xi_{i}^{-1}g) \cdot \sigma_{i}(u_{i}g) \cdot \sigma_{i}(u_{i}g)^{-1} \cdot (\xi_{i}g) \cdot (hg) \cdot (\xi_{i}^{-1}g) \cdot \sigma_{i}(u_{i}g)
\cdot \sigma_{i}(u_{i}g)^{-1} \cdot (\xi_{i}g) \cdot (\sigma_{i}(u)^{-1}g) = \\
= (\sigma_{i}(u)g) \cdot (\xi_{i}g) \cdot (\sigma_{i}(u)^{-1}g) = \psi_{i}(u,g) \cdot g.
\end{multline*}
\boom
The transition functions of these charts are $\alpha_{ij} \co P_{ij} \rightarrow \Aut(H)$ defined by
\[
\alpha_{ij}(u)(h) = s_{ij}(u) \cdot h \cdot s_{ij}(u)^{-1}.
\]
As far as the isometablicity of the respective transition functions is concerned, the following proposition is a straightforward calculation.
\begin{prop}
The transition functions $\alpha_{ij}$ are isometablic in the sense
\begin{eqnarray}\label{isomtransgpd}
\alpha_{ij}(ug)(\rho_{jj}(g^{-1})(h)) = \rho_{ii}(g^{-1})(\alpha_{ij}(u)(h)).
\end{eqnarray}
\end{prop}

\begin{thm}\label{classPBGLGB}
Let $P(M,G)$ be a principal bundle, ${\cal{P}} = \{ P_{i} \}_{i \in I}$ be an open cover of $P$ by principal bundle charts and $H$ a Lie group. Let ${\cal{\rho}} = \{ \rho_{i} \}_{i \in I}$ be a family of actions of $G$ on $H$. Given a cocycle ${\cal{\alpha}} = \{ \alpha_{ij} \co P_{ij} \rightarrow \Aut(H) \}_{i,j \in I}$ which satisfies (\ref{isomtransgpd}), there exists a PBG-Lie group bundle over $P(M,G)$ with transition functions the given ones.
\end{thm}

\pf Let $F_{i} = U_{i} \times H$ and on the union of the $F_{i}$ define the following equivalence relation:
\[
(i,(u_{1},h_{1})) \sim (j,(u_{2},h_{2})) \Leftrightarrow u_{1} = u_{2} = u~{\mbox and}~h_{2} = \alpha_{ij}(u)(h_{1}).
\]
This is an equivalence relation because we assumed that the $\alpha_{ij}$s form a cocycle. Denote the quotient set by $F$ and equivalence classes $\langle i,(u,h) \rangle$. Define a map $\pi \co F \rightarrow P$ by $\pi\langle i,(u,h) \rangle = u$ and a $G$-action by
\[
\langle i,(u,h) \rangle \cdot g = \langle i,(ug,\rho_{i}(g^{-1})(h)) \rangle.
\]
It is easy to see that the map $\psi_{i} \co P_{i} \times H \rightarrow \pi^{-1}(P_{i})$ defined by $(u,h) \mapsto \langle i,(u,h) \rangle$ is an equivariant bijection. Give $F$ the smooth structure induced from the manifolds $P_{i} \times H$ via the $\psi_{i}$s. Clearly $F \rightarrow P(M,G)$ is a PBG-Lie algebra bundle, and its transition functions are
\[
\psi_{i,u}^{-1} (\psi_{j,u}(h)) = \psi_{i,u}^{-1}(\langle j,(u,h) \rangle) = \psi_{i,u}^{-1}(\langle i,(u,\alpha_{ij}(u)(h)) \rangle) = \alpha_{ij}(u)(h).
\]
\boom

It will be shown in Section 6 that the construction of a PBG-LGB given in \ref{classPBGLGB} is well defined.

The family of $G$-actions $\{ \rho_{ij} \}_{i,j \in I}$ arises naturally from the local flat basic connections that every PBG-groupoid has. A remarkable result, which is presented here, is that these actions are really only local expressions of the $G$-action on the groupoid. We prove this for the subset of the $\rho_{ij}$s for which $i = j$. This is enough, as it was shown in (\ref{indepj}) that these actions determine the whole family. To this end, it is necessary to establish the notion of an action groupoid.
\begin{df}\label{actgpd}
Given a manifold $M$ together with a right action of a Lie group $G$ on $M$, the {\em action groupoid} $M \ract G \gpd M$ associated with this action is the product manifold $M \times G$, together with the following groupoid structure:
\begin{enumerate}
\item The source map is $\alpha(x,g) = x$, and the target map $\beta(x,g) = xg$.
\item Multiplication is defined by $(xg,h)\cdot (x,g) = (x,gh)$.
\item The unit element over any $x \in M$ is $1_{x} = (x,e_{G})$.
\item The inverse of an element $(x,g) \in M \ract G$ is $(xg,g^{-1})$.
\end{enumerate}
\end{df}
Note that the action groupoid is transitive if and only if the $G$-action on $M$ is transitive.

Now suppose given a PBG-groupoid $\Omega \gpd P(M,G)$ and a cover ${\cal{P}} = \{ P_{i} \}_{i \in I}$ of $P$ by principal bundle charts. For every $i \in I$, consider the action groupoid $P_{i} \ract G \gpd P_{i}(U_{i},G)$ and define a map $\tilde{\rho}_{i} \co P_{i} \ract G * I\Omega_{P_{i}} \rightarrow I\Omega_{P_{i}}$ by
\[
\Tilde{\rho}_{i}((u,g),\eta \in \Omega_{u}^{u}) = \psi_{i}(ug,\rho_{ii}(g^{-1})(\psi_{i,u}^{-1}(\eta))).
\]
Obviously, $\pi(\Tilde{\rho}_{i}((u,g),\eta)) = ug = \beta(u,g)$ and $\Tilde{\rho}_{i}((u,e_{G}),\eta) = \eta$. It is easily verified that
\[
\Tilde{\rho}_{i}((ug_{1},g_{2}) \cdot (u,g_{1}),\eta) = \Tilde{\rho}_{i}((ug_{1},g_{2}),\Tilde{\rho}_{i}(u,g_{1}),\eta).
\]
Also, each $\Tilde{\rho}_{i}(u,g)$ is an automorphism of $\Omega_{u}^{u}$, therefore it is a representation of the Lie groupoid $P_{i} \ract G$ on the Lie group bundle $I\Omega_{P_{i}}$, in the sense of \cite{Mac:class}. The following proposition allowes us to "glue" the $\Tilde{\rho}_{i}$s together to a global map.
\begin{prop}
For all $i,j \in I$ such that $P_{ij} \neq \emptyset$, $u \in P_{ij}$, $g \in G$ and
$\eta \in \Omega_{u}^{u}$ we have
\[
\Tilde{\rho}_{i}((u,g),\eta) = \tilde{\rho}_{j}((u,g),\eta).
\]
\end{prop}

\pf The isometablicity of the $\alpha_{ij}$'s gives:
\begin{multline*}
\Tilde{\rho}_{i}((u,g),\eta) =
\psi_{i}(ug,\rho_{ii}(g^{-1})(\psi_{i,u}^{-1}(\eta))) =
\psi_{i}(ug,\rho_{ii}(g^{-1})(\alpha_{ij}(u)(\psi_{j,u}^{-1}(\eta)))) = \\
= \psi_{i}(ug,\alpha_{ij}(ug)(\rho_{jj}(g^{-1})(\psi_{i,u}^{-1}(\eta)))) =
\psi_{j}(ug,\rho_{jj}(g^{-1})(\psi^{-1}_{i,u}(\eta))) 
= \tilde{\rho}_{j}((u,g),\eta).
\end{multline*}
\boom

Now we can define $\rho : (P \ract G) * I\Omega \rightarrow I\Omega$ by
$\rho((u,g),\eta \in \Omega_{u}^{u}) = \Tilde{\rho}_{i}((u,g),\eta)$, if $u \in P_{i}.$ The
previous proposition shows that it is well defined. More than that, it is a representation because each
$\tilde{\rho}_{i}$ is. As a matter of fact, $\rho$ is a lot simpler than it seems. Since the charts 
$\{ \psi_{i} \}_{i \in I}$ are isometablic we have:
\[
\rho((u,g),\eta) = \psi_{i}(ug,\rho_{ii}(g^{-1})(\psi_{i,u}^{-1}(\eta))) =
\psi_{i}(u,\psi_{i,u}^{-1}(\eta)) \cdot g = \eta \cdot g.
\]
So $\rho$ is, in fact, just the PBG structure of $I\Omega.$

Conversely, it is possible to retrieve the local representations $\{ \rho_{ii} \}_{i \in I}$ from the PBG structure of $I\Omega$. Suppose $\{ \sigma_{i} : P_{i} \rightarrow \Omega_{u_{0}} \}_{i \in I}$ is a family of sections of $\Omega$. Consider the charts $\psi_{i} : P_{i}
\times H \rightarrow I\Omega_{P_{i}}$ defined as $\psi_{i,u}(h) =
I_{\sigma_{i}(u)}(h)$ and define $\Tilde{\rho}_{i} : P_{i} \ract G \rightarrow Aut(H)$ 
by
\[
\Tilde{\rho}_{i}(u,g)(h) = \psi_{i,ug}^{-1}(\psi_{i,u}(h) \cdot g)
\]
for all $g \in G$, $h \in H$ and $u \in P_{i}$. This is a morphism of Lie groupoids over
$P_{i} \rightarrow \cdot$. For every $i \in I$ choose $u_{i} \in P_{i}$ and define
\[
\rho_{ii}(g^{-1})(h) = \Tilde{\rho}_{i}(u_{i},g)(h) = \psi^{-1}_{i,u_{i}g}
(\psi_{i,u}(h) \cdot g).
\]
Then,
\begin{multline*}
\rho_{ii}(g^{-1})(h) = I_{\sigma_{i}(u_{i}g)^{-1}}(I_{\sigma_{i}(u_{i})}(h)
\cdot g) = \sigma_{i}(u_{i}g)^{-1} \cdot (\sigma_{i}(u_{i})g) \cdot (hg) \cdot
(\sigma_{i}(u_{i})^{-1}g) \cdot \sigma_{i}(u_{i}g).
\end{multline*}
The latter is exactly the original definition of the $\rho_{ii}$'s. Since the $\rho_{ii}$s determine the $\rho_{ij}$s, the previous considerations are the proof of the following theorem:
\begin{thm}\label{localisom}
Given a PBG-groupoid $\Omega \gpd P(M,G)$, the representations $\{ \rho_{ij} \}_{i \in I}$ are local expressions of the PBG structure of $I\Omega$. 
\end{thm}

\section{The classification of transitive PBG-groupoids}

This section deals with a single result: It is shown that the isometablic transition functions classify transitive PBG--groupoids. 

\begin{thm}\label{classifPBGgpds}
Let $P(M,G)$ be a principal bundle and ${\cal{P}} = \{ P_{i} \}_{i \in I}$ an open cover of $P$ by
principal bundle charts. Consider a Lie group $H$ and a family of actions
${\cal{\rho}} = \{ \rho_{ij} \}_{i,j \in I}$ of $G$ on $H$ which has the property of the
cocycle morphism. Given a ${\cal{\rho}}$-isometablic cocycle
$\{ s_{ij} : P_{ij} \rightarrow H \}_{i,j \in I}$ there is a PBG-groupoid $\Omega$ over
$P(M,G)$ whose PBG-Lie group bundle $I\Omega$ of orbits is the one produced by
$\{ \alpha_{ij} = I_{s_{ij}}\}_{i,j \in I}.$
\end{thm}

\pf For every $i,j \in I$ consider the sets $\Sigma_{i}^{j} = P_{i} \times H \times P_{j}$ and
let $\Sigma = \bigcup_{i,j \in I} \Sigma_{i}^{j}$. Consider the equivalence relation
\begin{multline*}
(i,u,h,v,j) \sim (i',u',h',v',j') \Longleftrightarrow u = u',~v = v' ~{\mbox and}~
h' = s_{i'i}(u) \cdot h \cdot s_{jj'}(v).
\end{multline*}
Then it is shown in \cite[II 2.19]{LGLADG} that the following defines a groupoid structure on the
quotient $\Omega = \Sigma/\sim$: The source and target projections are
$\langle i,u,h,v,j \rangle = v,~\langle i,u,h,v,j \rangle = u$, the object inclusion map is
$1 : u \mapsto 1_{u} = \langle i,u,e_{H},u,i \rangle$ (any $i \in I$ such that $u \in P_{i}$),
and the multiplication is
\[
\langle i,u,h_{1},v,j_{1} \rangle \cdot \langle j_{2},v,h_{2},w,k \rangle =
\langle i,u,h_{1} \cdot s_{j_{1}j_{2}}(v) \cdot h_{2},w,k \rangle.
\]
The inversion is $\langle i,u,h,v,j \rangle^{-1} = \langle j,v,h^{-1},u,i \rangle$. This
groupoid becomes a PBG-groupoid with action:
\[
\langle i,u,h,v,j \rangle \cdot g = \langle i,ug,\rho_{ij}(g^{-1})(h),vg,j \rangle.
\]
This is well defined because if $\langle i,u,h,v,j \rangle = \langle i',u,h',v',j' \rangle$ then
$h' = s_{i'i}(u) \cdot h \cdot s_{jj'}(v)$. The cocycle morphism condition then
gives:
\begin{multline*}
\rho_{i'j'}(g^{-1})(h') = \rho_{i'i}(g^{-1})(s_{i'i}(u)) \cdot
\rho_{ij}(g^{-1})(h) \cdot \rho_{jj'}(g^{-1})(s_{jj'}(u)) = \\
= s_{i'i}(ug) \cdot \rho_{ij}(g^{-1})(h) \cdot s_{jj'}(ug).
\end{multline*}
So, $\langle i,u,h,v,j \rangle \cdot g = \langle i',u,h',v,j' \rangle \cdot g$. It is
straightforward that this action makes $\Omega$ a PBG-groupoid. For instance, we prove here
that this action preserves the multiplication. Again, because of the cocycle morphism property,
we have:
\begin{multline*}
(\langle i,u,h_{1},v,j_{1} \rangle \cdot \langle j_{2},v,h_{2},w,k \rangle) \cdot g =
\langle i,u,h_{1} \cdot s_{j_{1}j_{2}} \cdot h_{2},w,k \rangle \cdot g = \\
= \langle i,ug,\rho_{ik}(g^{-1})(h_{1} \cdot s_{j_{1}j_{2}} \cdot h_{2}),w,k
\rangle = \\
= \langle i,ug,\rho_{ij_{1}}(g^{-1})(h_{1}) \cdot
\rho_{j_{1}j_{2}}(g^{-1})(s_{j_{1}j_{2}}(v)) \cdot
\rho_{j_{2}k}(g^{-1})(h_{2}),wg,k \rangle = \\
= \langle i,ug,\rho_{ij_{1}}(g^{-1})(h_{1}) \cdot s_{j_{1}j_{2}}(vg) \cdot
\rho_{j_{2}k}(g^{-1})(h_{2}),wg,k \rangle = \\
(\langle i,u,h_{1},v,j_{1} \rangle \cdot g) \cdot (\langle j_{2},v,h_{2},w,k \rangle \cdot g).
\end{multline*}
\boom

\begin{prop} \label{welldefined}
Let $P(M,G)$ be a principal bundle, $\{ P_{i} \}_{i \in I}$ an open cover of $P$ by principal bundle charts, $H$ a Lie group
and ${\cal{\rho}'}, {\cal{\rho}}$ be two families of actions of $G$ on $H$ by cocycle morphisms
which are conjugate under a family of maps $r = \{ r_{i} : P_{i} \rightarrow H \}_{i \in I}$
such that $r_{i}(ug) = \rho_{ii}(g^{-1})(r_{i}(u)) \cdot r_{i}(u_{i}g)$ for all $u \in P_{i},g \in G$ and $i \in I$. Let $\{ s_{ij} \}_{i,j \in I}$ and $\{ {s'}_{ij} \}_{i,j \in I}$ be
${\cal{\rho'}}$-isometablic and ${\cal{\rho}}$-isometablic systems of transition data over
$\{ P_{i} \}_{i \in I}$ with values in $H$ respectively which are equivalent under the family
of maps $r$. Let $\Omega'$ and
$\Omega$ be the PBG-groupoids constructed from $\{ s_{ij} \}_{i,j \in I}$ and
$\{ {s'}_{ij} \}_{i,j \in I}$ respectively. Then the map $\phi : \Omega' \rightarrow \Omega$
defined by
\[
\langle i,u,v,h \rangle \mapsto \langle i,u,r_{i}(u) \cdot h \cdot r_{j}(v)^{-1},
v,j \rangle
\]
is an isomorphism of PBG-groupoids over $P(M,G)$.
\end{prop}

\pf~It is shown in \cite[II 2.19]{LGLADG} that $\phi$ is an isomorphism of Lie groupoids.
To show that it is an isomorphism of PBG-groupoids, take any $g \in G$. Then
\begin{multline*}
\phi(\langle i,u,h,v,j \rangle \cdot g) = \phi(\langle i,ug,{\rho'}_{ij}(g^{-1})(h),vg,j \rangle) = \\
= \langle i,ug, r_{i}(ug) \cdot r_{i}(u_{i}g)^{-1} \rho_{ij}(g^{-1})(h) \cdot r_{j}(u_{j}g) \cdot r_{j}(vg)^{-1} ,vg,j \rangle = \\
= \langle i, ug,  \rho_{ii}(g^{-1})(r_{i}(u)) \cdot r_{i}(u_{i}g) \cdot r_{i}(u_{i}g)^{-1}
\rho_{ij}(g^{-1})(h) \cdot r_{j}(u_{j}g) \cdot r_{j}(u_{j}g)^{-1} \rho_{jj}(g^{-1})(r_{j}(v)^{-1}),vg,j \rangle = \\
= \langle i,ug, \rho_{ij}(g^{-1})(r_{i}(u) \cdot h \cdot r_{j}(v)^{-1}), vg,j \rangle 
= \phi(\langle i,u,h,v,j \rangle) \cdot g.
\end{multline*}
\boom

\section{Isometablic transition data}
\begin{sloppypar}
Let us move to the Lie algebroid level for a while. In \cite[IV\S4]{LGLADG}, it is shown that a transitive Lie algebroid $L \inj A \surj TM$ is locally described by the following data: If $\hoh$ denotes the fibre type of $L$, then for a simple open cover $\{ U_{i} \}_{i \in I}$ of $M$, there exists a family of differential-2-forms $\chi = \{ \chi_{ij} \co TU_{ij} \times TU_{ij} \rightarrow U_{ij} \times \hoh \}_{i,j \in I}$and a cocycle $\alpha = \{ \alpha_{ij} \co U_{ij} \rightarrow \Aut(\hoh) \}_{i,j \in I}$ such that
\end{sloppypar}
\begin{enumerate}
\item The $\chi_{ij}$s are Maurer-Cartan forms, i.e. $\delta\chi_{ij} + [\chi_{ij},\chi_{ij}] = 0$, whenever $U_{ij} \neq \emptyset$,
\item $\chi_{ik} = \chi_{ij} + \alpha_{ij}(\chi_{jk})$, whenever $U_{ijk} \neq \emptyset$,
\item $\Delta(\alpha_{ij}) = \ad \circ \chi_{ij}$, whenever $U_{ij} \neq \emptyset$.
\end{enumerate}
The $\alpha_{ij}$s here are the transition functions of the Lie algebra bundle $L$. The notation $\Delta$ stands for the Darboux derivative. More than that, it is shown that this data classifies transitive Lie algebroids.

Since transitive Lie groupoids differentiate to transitive Lie algebroids, it is reasonable to expect that so does the respective classification data. Mackenzie in \cite[III\S5]{LGLADG}, gives a full account of this process, however it is expected that the transition functions that classify a transitive Lie groupoid can be reformulated in a fashion which makes their correspondence to the pair $(\chi,\alpha)$ on the algebroid level immediate. 

In this section we give this reformulation for transitive PBG-groupoids. For any PBG-groupoid $\Omega \rightarrow P(M,G)$ such that the fibre bundle of the associated PBG-Lie group bundle $I\Omega$ is $H$, we have the following definition:

\begin{df}
The Lie groupoid morphisms $\chi_{ij} \co P_{ij} \times P_{ij} \rightarrow H$ defined by
\[
\chi_{ij}(u,v) = s_{ij}(u) \cdot s_{ji}(v)
\]
(over the map $P_{ij} \rightarrow \cdot$) are called {\em transition morphisms}.
\end{df}
Let us see now how the transition morphisms intertwine with the transition functions $\alpha_{ij}$.
\begin{prop}
The transition morphisms $\chi_{ij}$ and the transition functions $\alpha_{ij}$ satisfy:
\begin{enumerate}
\item $\chi_{ik}(u,v) = \chi_{ij}(u,v) \cdot \alpha_{ij}(v)(\chi_{jk}(u,v))$
\item For a choice of $u_{ij} \in P_{ij}$, 
\[
\alpha_{ij}(u) = I_{\chi_{ij}(u,u_{ij})} \circ I_{s_{ij}(u_{ij}))}.
\]
\item $\rho_{ii}(g^{-1})(\chi_{ij}(u,v)) = \chi_{ij}(ug,vg)$.
\end{enumerate}
\end{prop}
Again, the proof is straightforward. Note that these conditions differentiate to the respective ones on the Lie algebroid level. 
\begin{df}
Let $P(M,G)$ be a principal bundle, ${\cal{P}} = \{ P_{i} \}_{i \in I}$ a cover of $P$ by principal bundle charts, $H$ a Lie group and ${\cal{\rho}} = \{ \rho_{i} \}_{i \in I}$ a family of $G$-actions on $H$. Let ${\cal{\chi}} = \{ \chi_{ij} \co P_{ij} \times P_{ij} \rightarrow H \}_{i,j \in I}$ be a family of Lie groupoid morphisms and ${\cal{\alpha}} = \{ \alpha_{ij} : P_{ij} \rightarrow \Aut(H) \}_{i,j \in I}$ a cocycle, such that
\begin{enumerate}
\item $\rho_{ii}(g^{-1})(\chi_{ij}(u,v)) = \chi_{ij}(ug,vg)$,
\item $\alpha_{ij}(ug)(\rho_{jj}(g^{-1})(h)) = \rho_{ii}(g^{-1})(\alpha_{ij}(u)(h))$,
\item $\chi_{ik}(u,v) = \chi_{ij}(u,v) \cdot \alpha_{ij}(v)(\chi_{jk}(u,v))$,
\item For a choice of $u_{ij} \in P_{ij}$, 
\[
\alpha_{ij}(u) = I_{\chi_{ij}(u,u_{ij})} \circ I_{s_{ij}(u_{ij}))}.
\]
\end{enumerate}
Then the pair $({\cal{\chi}},{\cal{\alpha}})$ is called a {\em ${\cal{\rho}}$-isometablic system of transition data over $P(M,G)$ with values in $H$}.
\end{df}
Let us now examine the relation of systems of transition data when we start with different families of flat isometablic basic connections. Denote $({\cal{\chi}},{\cal{\alpha}})$ and $({\cal{\chi'}},{\cal{\alpha'}})$ the respective systems of isometablic transition data. Again, the proof of the following proposition is just a matter of calculations.
\begin{prop}
Two ${\cal{\rho}}$-isometablic and ${\cal{\rho'}}$-isometablic systems of transition data $({\cal{\chi}},{\cal{\alpha}})$ and $({\cal{\chi'}},{\cal{\alpha'}})$ respectively are related by
\begin{eqnarray}\label{systrelchi}
{\chi'}_{ij}(u,v) = r_{i}(u)^{-1}[\chi_{ij}(u,v) \cdot \alpha_{ij}(v)(r_{i}(u) \cdot r_{j}(v)^{-1})] \cdot r_{i}(v)
\end{eqnarray}
and
\begin{eqnarray}\label{systrelalpha}
{\alpha'}_{ij}(u) = I_{r_{i}(u)^{-1}} \circ \alpha_{ij}(u) \circ I_{r_{j}(u)}
\end{eqnarray}
\end{prop}
\begin{df}
Two isometablic systems of transition data which satisfy ({\ref{systrelchi}}) and (\ref{systrelalpha}) are called {\em equivalent}.
\end{df}
Finally we prove that the PBG-Lie group bundles induced by equivalent transition functions are isomorphic, showing thus that the classification of PBG-Lie group bundles we gave in \ref{classPBGLGB} is well defined.
\begin{thm}
Let $P(M,G)$ be a principal bundle, ${\cal{P}} = \{ P_{i} \}_{i \in I}$ a cover of $P$ by principal bundle charts and $H$ a Lie group. Let ${\cal{\rho}} = \{ \rho_{i} \}_{i \in I}$ and ${\cal{\rho'}} = \{ {\rho'}_{i} \}_{i \in I}$ be two families of actions of $G$ on $H$ such that
\begin{enumerate}
\item $\rho_{i}(g^{-1})(h_{1}h_{2}) = \rho_{i}(g^{-1})(h_{1}) \cdot \rho_{i}(g^{-1})(h_{2})$
\item There exists a family of maps $\{ r_{i} \co P_{i} \rightarrow H \}_{i \in I}$ which are ${\cal{\rho}}$-isometablic (i.e. $r_{i}(ug) = \rho_{i}(g^{-1})(r_{i}(u)) \cdot r_{i}(u_{i}g)$) such that
\[
{\rho'}_{i}(g^{-1})(h) = r_{i}(u_{i}g)^{-1} \cdot \rho_{i}(g^{-1})(h) \cdot r_{i}(u_{i}g).
\]
\end{enumerate}
If ${\cal{\alpha}}$ and ${\cal{\alpha'}}$ are cocycles which satisfy (\ref{systrelalpha}) which give rise to the PBG-Lie group bundles $F$ and $F'$ respectively, then the map $\phi \co F \rightarrow F'$
\[
\langle i,(u,h) \rangle \mapsto \langle i,(u,r_{i}(u)^{-1} \cdot h \cdot r_{i}(u)) \rangle
\]
is an isomorphism of PBG-Lie algebra bundles.
\end{thm}
The proof of this is analogous to the one given in \ref{welldefined}.

\section{Isometablic cohomology}

In this section we give a formulation of the cohomology that classifies PBG-groupoids. In general, consider a principal bundle $P(M,G)$, a cover ${\cal{P}} = \{ P_{i} \}_{i \in I}$ of $P$ by principal bundle charts and a Lie group $H$. We also suppose given a family ${\cal{\rho}} = \{ \rho_{ij} \}_{i,j \in I}$ of $G$-actions on $H$ with the property of the cocycle morphism. 

For $n \geq 3$ we denote by $\check{C}^{n}_{G}(P,H)$ the set of differentiable maps $e_{i_{0},\ldots, i_{n}} \co P_{i_{0},\ldots,i_{n}} \rightarrow H$ such that for every $u \in P_{i_{0},\ldots,i_{n}}$ and $g \in G$ we have:
\begin{enumerate}
\item $e_{i_{0},\ldots, i_{n}}(ug) = \rho_{i_{n-1},i_{n-2}}(g^{-1})(e_{i_{0},\ldots, i_{n}}(u))$, if $n$ is odd and
\item $e_{i_{0},\ldots, i_{n}}(ug) = \rho_{i_{n-1},i_{n-3}}(g^{-1})(e_{i_{0},\ldots, i_{n}}(u))$, if $n$ is even.
\end{enumerate}
For $n = 0$ define $\check{C}^{0}_{G}(P,H)$ to be the set of $e_{i} \co P_{i} \rightarrow H$ such that $e_{i}(ug) = \rho_{ii}(g^{-1})(e_{i}(u))$. For $n = 1$ define $\check{C}^{1}_{G}(P,H)$ to be the set of $e_{ij} \co P_{ij} \rightarrow H$ such that $e_{ij}(ug) = \rho_{ij}(g^{-1})(e_{ij}(u))$. Finally, define $\check{C}^{2}_{G}(P,H)$ to be the set of $e_{ijk}$s such that $e_{ijk}(ug) = \rho_{jj}(g^{-1})(e_{ijk}(u))$ and identify $\check{C}^{-1}_{G}(P,H)$ with $H$.

Then the usual \v{C}ech differential $\delta \co \check{C}^{n}(P,H) \rightarrow \check{C}^{n+1}(P,H)$ defined by
\[
\delta(e)_{i_{0},\ldots,i_{n}} = \Pi_{k = 0}^{n} [e_{i_{0},\ldots,\Hat{i_{k}},\ldots,i_{n}}]^{(-1)^{k + 1}}
\]
is {\em isometablic} in the sense that
\begin{enumerate}
\item $\delta(e)_{i_{0},\ldots,i_{n}} (ug) = \rho_{i_{n-1}i_{n-2}}(g^{-1})
(\delta(e)_{i_{0},\ldots,i_{n}}(u))$ if $n$ is odd and 
\item $\delta(e)_{i_{0},\ldots,i_{n}} (ug) = \rho_{i_{n-1}i_{n-3}}(g^{-1})
(\delta(e)_{i_{0},\ldots,i_{n}}(u))$ if $n$ is even.
\end{enumerate}
\begin{df}
The cohomology of the complex 
\[
... \stackrel{\delta}{\rightarrow} \check{C}^{n}_{G}(P,H) \stackrel{\delta}{\rightarrow} \check{C}^{n+1}_{G}(P,H) \stackrel{\delta}{\rightarrow} ... 
\]
is called {\em isometablic \v{C}ech} cohomology and denoted by $\check{H}^{n}_{G}(P,H)$.
\end{df}

The next theorem follows immediately from \ref{classifPBGgpds}.

\begin{thm}
With the notation above, PBG-groupoids are classified by $\check{H}^{1}_{G}(P,H)$.
\end{thm}

\bibliographystyle{plain}

Current address:
\\ 
{\tt 
Iakovos Androulidakis
\\ 
19 Drosi Street
\\ 
11474 Athens
\\ 
Greece
}
\end{document}